\input amssym.tex
\input cyracc.def

\magnification =1200
\topskip = .5truecm
\vsize = 210true mm
\hsize =153 true mm
\baselineskip=15truept
\voffset=2\baselineskip
\overfullrule = 0pt
\parindent=0pt

\font\bfxx=cmbx12 at 15pt
\font\bfx=cmbx12

\font\ninesc=cmcsc9
\font\cmsc=cmcsc10

\newfam\cyrfam
\font\fivecyr=wncyr5%
\font\sevencyr=wncyr7%
\font\tencyr=wncyr10%

\def\cyrf{\fam\cyrfam\sevencyr\cyracc}

\textfont\cyrfam=\tencyr 
\scriptfont\cyrfam=\sevencyr%
\scriptscriptfont\cyrfam=\fivecyr

\font\bfn=cmbx9

\def\inv#1{#1^{-1}}
\def\recip#1{\frac{1}{#1}}
\def\frac#1#2{{#1\over #2}}

\def\half{\recip{2}}

\def\intdl#1#2#3#4{\int_{#1}^{#2}#3\,{\rm d}#4}

\def\ref#1{[{\sl #1\/}]}

\newcount\notenumber

\def\noter#1{\advance\notenumber by 1
$^{\the\notenumber}$\footnote{}{$^{\the\notenumber}$\ #1}}
\outer\def\subsectionsc#1\par{\hfill\smallskip\message{
#1}\leftline{\ninesc #1}\nobreak\smallskip\noindent } %

\outer\def\proclaimsc#1#2{\par\medbreak\noindent{\ninesc #1 \enspace}{\sl
#2}\par\ifdim\lastskip<\medskipamount \removelastskip\penalty55\medskip\fi}%

%
%

\newcount\sectnumber%
\def\sectn#1 {\advance \sectnumber by 1{\par\bigskip\noindent{\bfx
\the\sectnumber\ #1}}\par}%

\newcount\subsnumber%
\def\subsn#1 {\advance \subsnumber by 1{\par\medskip\noindent{\cmsc\the\sectnumber.\the\subsnumber\ 
#1}\quad}}%
\outer\def\subsectionsc#1\par{\hfill\smallskip\message{
#1}\leftline{\ninesc #1}\nobreak\smallskip\noindent } %

\newcount\subsanumber%
\def\asubsn#1 {\advance \subsanumber by 1{\par\medskip\noindent{\cmsc\the\sectnumber.\the\subsanumber\ 
#1}\quad}}%

\newcount\subsbnumber%
\def\bsubsn#1 {\advance \subsbnumber by 1{\par\medskip\noindent{\cmsc\the\sectnumber.\the\subsbnumber\ 
#1}\quad}}%

\def\napa#1{\hfill\break\noindent{\bfn #1} }%
\pageno = 1

\null
\vskip -1truecm

\centerline{\bfxx More About Trigonometric	Series and  Integration	}

\bigskip

\sectn{ Introduction}

In a recent interesting paper by  Gluchoff,  \ref{G}, the development of the  theory of integration	 was
related to the needs of trigonometric	series. The paper ended with the
introduction of the Lebesgue	 integral	at the beginning of the century. As the Lebesgue	
integral	is the integral	of  everyday mathematics this was a natural place to
stop; the needs of trigonometric series had produced a tool of first-rate importance
for the whole of mathematical analysis. 

However there is a question left over from [G]; it  is the problem of determining the
coefficients of an everywhere convergent trigonometric series 
in terms of its sum, the so called {\sl coefficient problem for convergent
trigonometric	series\/}, that we will call shortly, {\sl the coefficient problem\/}. Although the need for
further developments in the theory of integration	 had disappeared, one of the basic problems of trigonometric	
series was not solved  by the Lebesgue	 integral.  As a  result,  right up to the present day  further integrals,
 usually called {\sl\ trigonometric	 integrals}, have been  introduced in order to solve  the coefficient problem. These
integrals  are
 more general than the Lebesgue integral.

 Although such  integrals have   not found other uses it might be of
some interest to continue the story started in \ref{G} by describing some of this work.
\smallskip
As this article  continues the discussion started in [G] reference will be made to the results and references there.  
\smallskip
A knowledge of integrals more general than that of
 Lebesgue is important; see for instance \ref{Bu2; Go; Ho; L; P}. The first reference was written to put this
information, some of which is rather technical, in an easily available form. It is almost essential reading for the
full understanding of the present topic, and use will be made of the results and references quoted there. 

Finally certain terms that are defined in the  the Appendix, Section 11, will be written as {\tt smooth}.
\bigskip
\sectn{The Coefficient Problem}

The  coefficient problem for convergent
trigonometric	series, or just  the coefficient problem, is the following.
\proclaimsc{}{ If for all $x, -\pi	\le x\le \pi	$,  the series 
$$
\frac{a_0}{2}\, +\, \sum_{k=1}^{\infty} a_k \cos kx + b_k\sin kx,\eqno(T)
$$
converges, to $f(x)$ say, 
  then this representation of  the function	$f$  is unique.\hfil\break
This implies that the coefficients of the series, $a_0,a_1,\ldots, b_1,
\ldots$, are completely determined by the sum  function, $f$.\hfil\break
Calculate these coefficients  from this function.}\par
 {\sl }
\smallskip
The fact that only one series of the form (T) can converge to the function $f$ is a famous result
due to Cantor; see \ref{G I, IV}, where references  are given for proofs.

As explained in \ref{G}  if the series (T) converges
uniformly then an argument that goes back to Fourier, and which is included in
all first  courses on Fourier analysis,  shows that 
$$
\eqalign
{
 a_k = &\frac{1}{\pi}\intdl{-\pi}{\pi}{f(x)\cos kx}{x},\; k=0,1,\ldots,\cr
 b_k  = &\frac{1}{\pi}\intdl{-\pi}{\pi}{f(x)\sin kx}{x},\; k=1,\ldots;\cr
}
\eqno(F)
$$
see \ref{G(1)}.  

In fact the argument that gives (F) only requires that the integrals
in (F) have a meaning, and term-by-term  integration
of the series (T) is justified; that is the integral of $f$  is obtained by integrating the series (T);  see \ref{G, pp.7--8}. In the
case of uniform convergence of a series of continuous functions it is well known that a series can be integrated term-by term.

Uniform convergence of the series  (T) implies that the sum $f$ is continuous  and so  the
integrands in (F) are   Riemann, or ${\cal R}$-,  integrable. It is then  natural to ask two questions. 
 \smallskip
{\sl
\item\item{(I)} Do the formul\ae\  (F) hold under the assumption that the  sum  $f$ of the series
(T) is ${\cal R}$-integrable?
\item\item{(II)} Is any $f$ that is the sum of a series (T) ${\cal R}$-integrable?
\/}
\smallskip
The answer to (I) is yes:  if the series  (T) converges everywhere to an ${\cal R}$-integrable
function then the coefficients of (T) are given by (F). This result is due to Du Bois-Reymond,\ref{Ho, vol.II, p.659}.

The answer to (II) is no: this was pointed out in \ref{G IV, p.18}. Lebesgue has shown that the sum a
function of a series  (T) need not be ${\cal R}$-integrable;  \ref{G, ref.15}.

\sectn{Why The Lebesgue Integral Will Not Do}  

The
next  step is to repeat the above questions (I) and (II)  with the  Lebesgue, or ${\cal L}$-,	 integral  instead of
the ${\cal R}$-integral.

 Again the answer  to the  question (I)  is is yes; if the series  (T) converges everywhere to a
${\cal L}$-integrable function then the coefficients of (T) are given by (F), where
now the integrals are of course Lebesgue integrals. This is  another result of Lebesgue; see 
\ref{G IV, pp.16-17; G, ref.16}.

Unfortunately the answer to the  question (II) is again no. We will show   why this is so, and at
the same time see just how bad  sum functions  of trigonometric series can be.

A natural candidate for the integral of $f$ is the sum of 
series you get on formally  integrating (T)  term-by-term;
$$
 	\Phi(x) = \frac{a_0}{2}x\, +\, \sum_{k=1}^{\infty} \frac{a_k \sin kx -
b_k\cos kx}{k}.\eqno(1)
$$
 This series is seemingly better behaved than the original series (T), since the
coefficients are essentially those of (T) divided by the terms of the sequence $\{1/n,\,
n=1,2,\ldots\}$. However this is  is not true in general.

 Consider the following example due to Fatou, where
$a_k= 0,\, k=0,1,2,\ldots$ and $b_k = 1/\log(k+1),\, k=1,2,\ldots$. Then (T) and (1) become

$$
\eqalignno
{
 &\sum_{k=1}^{\infty}\frac{\sin kx}{\log (k+1)},&(2)\cr
  \Xi(x)
=&-\sum_{k=1}^{\infty}\frac{\cos kx}{k\log (k+1)},&(3)\cr
}
$$
respectively.

 A simple application of  the Dirichlet convergence test,  \ref{K, pp.315, 316 Examples and
Applications 5},  shows that  the series (2)   converges for all $x$, to $\xi(x)$ say.

However the once integrated series in  (3)  does not converge when $x=0$. In other
words $\Xi$, the ``integral'' of $\xi$, is not defined at the origin. In particular it is not
continuous there.

 However all standard indefinite integrals are defined  everywhere, and are
continuous.  So  we must face the possibility that $\xi$ is not integrable in any standard sense. In
particular it cannot be ${\cal L}$-integrable.

Of course the above argument is only heuristic, as
(T) cannot necessarily be integrated term-by-term. However it can be made precise, and we will  give a sketch of the
proof; for details see \ref{D, vol.I, pp.42--43}.

 The  series  in (2), which as we saw
is convergent for all
$x$,   is  uniformly convergent on any
$[2k\pi +\epsilon, (2k+1)	\pi- \epsilon]$, where $0<\epsilon<\pi$, and $ k=0, \pm 1, \pm 2,
\ldots$;  see \ref{K, pp. 347, 349 Examples and Applications 4}.

The series obtained by the formally integrating (2), integrating
term-by-term, is   (3), which by a use of the  Dirichlet test 
 diverges for all $x= 2k\pi, k=0, \pm 1,\pm 2,\ldots$. Since (2) converges uniformly
  on  any $[\alpha, \beta]$ that does not contain
one of these points, that is,  on  any $[\alpha, \beta]$ that lies in some $[2k\pi +\epsilon, (2k+1)	\pi- \epsilon]$,  we get using
term-by-term integration that
$$
 \Xi(\beta) -\Xi(\alpha) = \int_{\alpha}^{\beta}\xi	,
$$
 where  the integral is taken in the Riemann sense. 

Suppose   $\xi$ were  ${\cal R}$-integrable and that $0< \alpha< \beta<\pi$,   then

$$
\int_{0}^{\beta}\xi =  \lim_{\alpha\to o}\int_{\alpha}^{\beta}\xi = \Xi(\beta) -\lim_{\alpha\to 0}\Xi(\alpha).
$$

However it can be shown that  for some  function    $m(x)$,
with
$xm(x)$ bounded away from zero,
 $$
\lim_{x\to 0}\biggl(\Xi(x)+ \sum_{k=1}^{m(x)}\frac{1}{k\log(k+1)}\biggr) = 0.
$$
So,  noting again that  $\displaystyle\sum\frac{1}{k\log(k+1)}$ diverges,  $\lim_{\alpha\to 0}\Xi(\alpha)$ does  not  exist .
\smallskip

Hence $\xi$ is not ${\cal R}$-integrable.  

It cannot be ${\cal L}$-integrable either, since  the above argument could be repeated using the Lebesgue integral.

More generally  $\xi$ cannot  be integrable in any
sense that extends the ${\cal R}$-integral, and has a continuous indefinite integral; see, for instance, the ${\cal D}^*$-integral in the next section.
\smallskip
The above argument applies to any series 
$$
\sum_{k=1}^{\infty} b_k\sin kx,
$$
 subject to the conditions
$$
b_1> b_2>\cdots;\qquad \lim_{k\to\infty}b_k = 0;\qquad
\sum_{k=1}^{\infty}\frac{b_k}{k} = \infty;
$$
as the references to \ref{K} show; see \ref{D, vol.I, pp. 42--43; Z, vol.I, pp.185--186}.
\sectn{ Why No Classical Integral Will Do}  

There are various extensions of the  Lebesgue integral and simplest of
these extensions  goes under various names that indicate the method of approach;  the narrow Denjoy, or
${\cal D}^*$-, integral, the  Perron, or ${\cal P}$-, integral and the Henstock-Kurzweil, or ${\cal HK}$-, 
integral; when we do not care about the
actual form we will refer to this integral as the ${\cal D}^*$-integral. This integral
extends the  ${\cal L}$-integral, in  that if a function is  ${\cal L}$-integrable then it is ${\cal D}^*$-integrable  and the
integrals are equal. In addition if a  function is bounded, or even just bounded below, and  is  ${\cal
D}^*$-integrable then  it is  ${\cal L}$-integrable. Finally the ${\cal D}^*$-integral, unlike the  ${\cal
L}$-integral, is not an {\sl absolute integral\/}, it is a {\sl non-absolute integral\/}; that is if $f$ is  ${\cal
D}^*$-integrable it does  not follow that
$|f|$ is. Finally  the ${\cal D}^*$-integral solves  the classical primitive problem, see Section 5 and
\ref{Bu2},  and the indefinite ${\cal D}^*$-integral, the ${\cal D}^*$-primitive, is continuous. The various approaches to this integral are
described in \ref{B2} and these  techniques will be used to define various trigonometric integrals .

From the last remarks, and the comments in the previous section, the sum of  (2) is not ${\cal D}^*$-integrable. So this integral, that is more
general than the ${\cal L}$-integral, does not solve the coefficient problem either; that is the answer to question (II) for 	 of Section 2 for this
integral  is again no.

However if the integral in the question (I) of Section 2 is
taken to be this more general  ${\cal D}^*$-integral the answer is yes; that is if the sum of
(T) is ${\cal D}^*$-integrable, the coefficients are given by (F), where the integral is if course the
${\cal D}^*$-integral; this result is due to Nalli, see [N;
Z, vol.II, pp. 83--86].

In \ref{Bu2, 6.3} a definition of an integral that is slightly more general than the ${\cal D}^*$-integral was given, the
${\cal DH}$-integral\noter{\sevenrm  Denjoy-Hin\v cin; Hin\v cin, {\cyrf Hinchin},  is also transliterated as
Khintchine; as a result this integral is sometimes called the  DK-integral.}. In this case the primitive is
differentiable almost everywhere and is $ACG$, while the  ${\cal D}^*$-primitive is differentiable almost everywhere
and is $ ACG^*$; these terms are defined  in \ref{Bu2, 6.3}.

Sklyarenko, \ref{S2;  Th3,  p.401}, has given an example of a ${\cal DH}$-integrable function, $f$ say, that is the sum of an everywhere
convergent series (T), but for which the coefficients are not given  by (F),  using the ${\cal DH}$-integral. In
particular 
$$
a_0\ne \recip{\pi}{\cal DH}\mkern-4mu-\mkern-8mu\int_{-\pi}^{\pi}\!f.
$$ 
 So with this slight extension of the ${\cal D}^*$-integral definition, the answer to question (I) becomes no. Of course the answer to question (II) is
also no, as the  ${\cal DH}$-primitive is continuous.
\smallskip
 However, it is
correct to say that
$\Phi$, the sum  of the series in (1) and often called the {\sl  Lebesgue  function of (T)\/}, is  the
natural primitive of $f$, the sum of the series (T). 

This
suggests that in order to have (F) hold in general a rather strange integral will be needed; one
whose primitive need not be continuous, or may even fail to be defined at some points.

The possibility of having an integration process with a primitive that is not
continuous was first considered by W.H.Young, \ref{Y}. It would appear that his
discussion was  motivated by theoretical considerations only, but it is
not completely impossible, given the important contributions that Young
made to the theory trigonometric series, that he was thinking of the need for such integrals
 in the coefficient problem.

\sectn{An Approach To The Coefficient Problem}

As  the above discussion  shows the coefficient problem can be solved if an  extra assumption
is made---uniform convergence of
the series  (T), the ${\cal R}$-, ${\cal  L}$-, or ${\cal D}^*$-integrability of the sum, $f$,   of the series (T). Such an assumption enables the
problem  to be solved but is unnecessary, since, as we have seen, the series (T)  when it converges need not be uniformly convergent, and its sum
function need not be  integrable in any of these senses..

What is needed is an integral, call it using a conceit due to Hardy, [H, p.7], the
{\sl Pickwickian integral\/}, or shortly the ${\cal  P}$i-{\sl integral\/}, for which  the answers to both
of the questions (I) and (II)  will be yes. That is:
\smallskip
{\sl 
\item\item{(A)}  all sum functions of everywhere convergent series (T)  will be ${\cal  P}$i-integrable;
\item\item{(B)} the coefficients of the series (T) are given by (F), where the integrals are ${\cal  P}$i-
integrals.
\/}
\smallskip
 Given a function $f$ for which the integrals in (F) exist in some sense, the
 series in  (T) can be written down. It is called the {\sl Fourier
series of
$f$\/}, whether the series converges to $f$ or not. The discussion of  two questions (I), and (II)  
above shows that if the series in  (T) converges to a ${\cal R}$-, ${\cal L}$-, or ${\cal
D}^*$integrable function $f$ then this series is the
Fourier series of $f$, more precisely it is a Fourier-Riemann series, a Fourier-Lebesgue
series, or a  Fourier-${\cal D}^*$ series of $f$.  The series (T) of the  example above of Sklyarenko is not  a Fourier-${\cal
DH}$ series; the sum function is  ${\cal DH}$- integrable but the coefficients calculated by (F) using this integral are not the coefficients of (T).
\smallskip
If we can define the ${\cal  P}$i- integral  having  the properties (A) and (B),
then (T) is  the Fourier-Pickwickian series of its sum, and the coefficient problem will have
been solved.

As we have noted the Pickwickian integral will be unusual as its primitive will not be
continuous, may even fail to exist at some  points.  The example of Fatou given above that has one bad point can, by the usual classical
technique of condensation of singularities, \ref{Ho, vol.II,  pp. 389--399},  be generalized to having an uncountable zero measure set  of 
discontinuities,  or where it is not defined; \ref{Bu2; D, vol.I, pp. 42; vol.IVb, pp.497--503}.

  However it cannot  get worse than
this. Consider the following facts.

 If (T) converges everywhere then  $a_k,
b_k\to 0$, see \ref{G IV, p.18}, and so by a simple application of the comparison test $\sum(a_k^2 + b_k^2)/k^2$  converges.
Hence by the Riesz-Fischer Theorem, see \ref{Z, vol. I, pp.127, 321},  if $\Phi$ is given by (1) then  $\Phi^2$ is ${\cal  L}$-
integrable; and this implies
$\Phi$ is finite almost everywhere; equivalently the series  (1)  converges almost
everywhere, that is
$\Phi$ can only fail to be defined on a set of measure zero.

Although the above remarks concentrate on the pathology of $\Phi$ they also indicate a way to
approach our problem. 

Suppose that we could obtain $\Phi$ from $f$, only knowing
 $f$  to be the sum of  (T).  Then by the above remarks $\Phi$ is ${\cal L}$-integrable, and  the series in  (1) is its 
Fourier-${\cal L}$ series. Hence the coefficients in (1)  are  given by (F), with $\Phi$ replacing $f$, the integrals being in the  Lebesgue sense
of course.  A comparison of (T) and (1) shows that their coefficients are simply related, so having the coefficients in (1) we have those of (T).

Since  $\Phi$ can  be considered as  an integral of $f$, we can also consider $f$ as a derivative  of
$\Phi$, in the Pickwickian sense say: $ D_{{\cal P}i}\Phi = f$.  Then the problem of obtaining $\Phi$ from
$f$ becomes the problem of  (a) defining what we mean by the ${\cal  P}$i-derivative, (b) showing that  $ D_{{\cal
P}i}\Phi = f$ and finally, (c) obtaining a method for inverting the  ${\cal P}$i-derivative. 

 In other words the coefficient problem can
be viewed as a generalization of the {\sl classical primitive problem\/}.

 The classical primitive  problem  is as follows.
\proclaimsc{}{If  on  $ [a,b],\, G'= g, G(a) = 0$ find  $G$ ,   or if $  G'= g$,  calculate $G(b) - G(a)$, from $g$.
}\par
A full discussion of this problem, together with the various solutions that have been given can be found in several
places; see in particular \ref{Bu2} where  there are further references. The various solutions, ${\cal D}^*$-
integral, the  Perron, or ${\cal P}$-, integral and the Henstock-Kurzweil, or ${\cal HK}$-, 
integral  are all
equivalent. The names indicate different approaches to the definition of the integral. As all  approaches---totalization for the ${\cal D}^*$-
integral, \ref{Bu2, 6.1},  major and minor functions for the ${\cal P}$- integral, \ref{Bu2, 6.2},  and Riemann sums for the  ${\cal HK}$-,
integral,  \ref{Bu2, 6.4}, --- are used in the coefficient problem it is important to know about these methods. In addition two very important
classes of functions are introduced when discussing these integrals,  the classes of continuous $ACG$ and continuous $ACG^*$ functions; \ref{Bu2,
6.3}.

 The various solutions of the coefficient problem that have been proposed are
obtained by replacing the ordinary derivative of the classical primitive  problem  by various
Pickwickian derivatives and then solving the  primitive problem for these derivatives.

It is clear that if we have inverted the Pickwickian derivative with a Pickwickian integral, that is 
$$
 D_{Pi}G = g\; \Longleftrightarrow\;G =  Pi\mkern-4mu-\mkern-10mu\int  g;
$$
 and if we also have, for $f$ the sum of (T), and $ \Phi$ given by (1),  that $ D_{Pi}\Phi = f$, then

$$
  Pi\mkern-4mu-\mkern-10mu\int_{\alpha}^{\alpha+2\pi}\!\!  f=\Phi(\alpha+2\pi) - \Phi(\alpha) =
\frac{a_0}{2}x\bigg|_{\alpha}^{\alpha+2\pi} = a_0\pi,
$$
 giving the first formula in (F),
$$
a_0 = \frac{1}{\pi}Pi\mkern-4mu-\mkern-10mu\int_{\alpha}^{\alpha+2\pi}  f.\eqno(4)
$$
Since  $\Phi$ may only be defined almost everywhere,  $\alpha$ is chosen so that the limits of integration are
points where $\Phi$ is defined.

 Once we have defined the correct ${\cal  P}$i- integral of course (4) will hold and so if $f$ is Sklyarenko's function
it will be integrable in both senses and
$$
{\cal DH}\mkern-4mu-\mkern-8mu\int_{-\pi}^{\pi}\!f\ne{\cal 
P}i\mkern-4mu-\mkern-8mu\int_{-\pi}^{\pi}\!f.
$$
 This gives  another indication of how awkward the ${\cal P}$i-integral has to be.  The class of ${\cal P}$i-integrable functions includes, by
Nalli's result quoted  section 4, all  ${\cal D}^*$-integrable function. So in  the class of ${\cal P}$i-primitives  includes  all continuous $ACG^*$
functions,  see \ref{Bu, 6.3}, but not, by Sklyarenko's example, all  continuous
$ACG$ functions. In fact the sum of the  Sklyarenko series has a ${\cal DH}$-primitive that is not its   ${\cal P}$i-primitive,  so the  ${\cal
P}$i-primitive is not an $ ACG $ function.

This example also shows that the difficulties  the  ${\cal P}$i-integral presents  are inherent in the 
coefficient problem.

Given that the ${\cal P}$i-integral calculates $a_0$, how do we calculate the rest of the coefficients of (T)?

One method is based on the theory of formal multiplication of trigonometric	series; \ref{Z, vol.I, pp.335--344; Ho. vol.II,
pp.585--587}. In this method the functions $ f(x)\sin nx, f(x)\cos nx$ are shown to  be the sums of
trigonometric	 series with constant terms $a_n/2, b_n/2$ respectively. Then the above
particular case (4) applies. 

Another method is to prove an integration by parts formula for the ${\cal P}$i-integral. This formula then plays
two roles; firstly it ensures that the functions $ f(x)\sin nx, f(x)\cos nx$ are ${\cal P}$i- integrable, 
and secondly it gives the required formul\ae. For instance:
$$
\eqalign
{
 Pi\mkern-4mu-\mkern-10mu\intdl{\alpha}{\alpha+2\pi}{\!f(x)\sin nx}{x} = &\Phi(x)\sin nx\big|_{\alpha}^{\alpha+\pi}-
\intdl{\alpha}{\alpha+2\pi}{\Phi(x)n\cos nx}{x}\cr
=& \pi b_n,\cr
}
$$
where the right a Lebesgue integral.

Both methods are useful as in some cases an integration by parts formula is not
easily obtained.

\sectn{The ${\cal P}${\rm i}-derivative Will Be  Symmetric Derivative}

 To see what role symmetry plays in the coefficient problem  consider the formula,  given in  \ref{G II, p.10}, for the
partial sums of the series in  (T) regarded as  the Fourier series of  $f$:

$$
s_n(x) \, =\, \frac{1}{\pi}\intdl{-\pi/2}{\pi/2}{f(x+ 2t)
\frac{\sin(2n+1)t}{\sin t}}{t}.\eqno(5)
$$
A simple change of variable in the integral (5) results in the alternative form
$$
s_n(x) \, =\, \frac{2}{\pi}\intdl{0}{\pi/2}{\biggl(\frac{f(x+
2t)+f(x-2t)}{2}\biggr) \frac{\sin(2n+1)t}{\sin t}}{t}.\eqno(6)
$$
The implication of (6) is that these partial sums  evaluated at $x$ depend
not on the value of $f$ at $x$ but on $f_e^x(t)$, the {\tt even part} of $f$ at $x$. In fact the heuristic argument
in
\ref{G II, pp.10--11} can be used to show that the series in  (T) converges to $\lim_{t\to 0}f_e^x(t) =(f(x+) +
f(x-))/2$ under the assumption that
$f$ is monotonic---the {\sl Dirichlet-Jordan test} for Fourier series; [Z, vol. I, p.57].

It should be remarked that on integrating a function, the even, respectively odd,  part of that function yields the odd,
respectively even part of its primitive.
	
It is natural to look at the Pickwickian derivatives  that use   the odd part of a function; that is, for example, a
derivative of  $\Phi$ of the kind defined  in (19).

 For a more detailed
discussion of the role of symmetry see [Th1,2], and the very important book [Th3].

\sectn{The Integrals of Denjoy and of James }

\subsn{Denjoy's $(T_{2,s})_o$-totalization.} The first solution of the coefficient problem, given by Denjoy, 
avoided the use of the primitive $\Phi$ because of the difficulties expected for the Pickwickian integral that have been
discussed above. Denjoy wanted an integral  with a continuous primitive  defined everywhere.

Denjoy  published
his result in a series of very short notes around 1920, but a full
version was only  given in the fifties as a result of a series of
lectures that Denjoy gave at Harvard on the invitation of McShane, [D].

The basic idea behind Denjoy's solution is to by-pass the first integral, $\Phi$,  and go
to the second integral
$$
 F(x) = \frac{a_0x^2}{4}  - \sum_{k=1}^{\infty}\frac{ a_k \cos kx +
b_k\sin kx}{k^2};\eqno(7)
$$
that is the series obtained by formally integrating  (T) twice.

As we remarked in Section 5 the coefficients of an everywhere convergent series  (T)  tend to zero;  this implies, by
a simple application of the comparison test,  that the series in (7) is uniformly, and absolutely,
convergent. Its sum  $F$ is then continuous, and   the series in (7)  is  the 
Fourier-${\cal R}$ series of
$F$. Hence the coefficients in this series, that are very simply related to those in
(T),  are given by (F), using the Riemann integral, and  with suitable changes in the integrand, using
$F$ instead of $f$. 

This would then calculate the coefficients of (T), except that we do not yet 
know how
$F$ is related to
$f$ except in this formal way---it is a sort of second integral.  Equivalently $f$ is a sort of second order derivative of $F$. It
follows from our discussion in Section 6 that we should look at  a derivative of the even part of $F$, and a careful
reading of the discussion in  \ref{G IV} will show in fact that the coefficient problem is
equivalent to a   
  primitive problem in which the Pickwickian derivative is  the {\tt Schwarz
derivative}; see (20).

 There is a famous
result of Riemann  discussed in \ref{G IV} that tells us how $f$, the sum of (T) is related to $F$, the sum in (7): namely, 
$D_s^2 F =f$.  A proof can be found in \ref{Z, vol.I, p.319}.

In our situation we know $f$ and wish to find $F$, the function of which it is the
second symmetric derivative. This is a  the generalization of the classical primitive problem used by Denjoy in his
solution of the coefficient problem.

It is known that if $ g= D_s^2G$ on $[a,b]$  for some continuous function $G$, then $G$ is unique
up to a linear function; this follows from the various results on trigonometric series given in
\ref{G}, but a direct proof is preferable and was given by Schwarz;  see \ref{Th3, p.12}. 

Hence if $G$ is any  continuous function such that $ 
D_s^2G= g$,   the unique  function that has the same property and  is zero at both
$a$ and at
$b$ is
$$
 \tilde G(x) = G(x)-G(a)-\frac{x-a}{b-a}\left(G(b) -G(a)\right)=(x-a)(x-b)\bigl[a,x,b;G\bigr],
$$
 where $\bigl[a,x,b; G\bigr]$ is the  {\tt second divided difference }of $G$ at $a,x,b$.

  In trying to solve the coefficient problem 
$F$, the sum of the series in of (7), is then regarded as  a second primitive of
$f$, in an integration process that inverts the Schwarz derivative.   

The process that Denjoy used was a extremely complicated extension of his totalization
process, see \ref{Bu2,  6.1} where other references are given. It defined what he called the  $({\cal T}_{2,s})_o$-integral and we write
$$
{\cal T}_{{(2,s)}_0}\mkern-4mu-\mkern-8mu\int_{(a,b)}^x g = \tilde G(x)=(x-a)(x-b)\bigl[a,x,b; G\bigr];\eqno(8).
$$
 the various subscripts 
are readily explained---$2$ for second order,
$s$ for symmetric, and $o$ for ordinary. This last needs bit of explanation; in Section 9 we will introduce an approximate derivative, here
however we are using ordinary derivative.

 Formula (8)  can be compared with 
$$
\int_a^x\!g = (x-a)[a,x; \Gamma]
$$
 that gives the  first  primitive
in terms of the {\tt first divided difference}.

If  $g$ is Lebesgue integrable  then   and so:
$$
{\cal T}_{{(2,s)}_0}\mkern-4mu-\mkern-8mu\int_{(a,b)}^x g = \frac{1}{b-a}\biggl\{(x-b)\intdl{a}{x}{(t-a)g(t)}{t}\, + \,
(x-a)\intdl{x}{b}{(t-b)g(t)}{t}\biggr\}\eqno(9).
$$
Formula (9) is easily shown if we note that
$$
G(x) - G(a) =\intdl{a}{x}{\intdl{a}{t}{f(u)}{u}}{t} =( x-a)\Gamma(x)-\intdl{a}{x}{(t-a)g(t)}{t},
$$
and another similar formula for $G(b)-G(x)$.

While the sum  $f$ of the series  (T) is, as we have seen, $ {\cal T}_{{(2,s)}_0}$-integrable on $[- \pi\pi]$, it is   not
necessarily ${\cal L}$-integrable. However, as in the case of the classical primitive problem, 
$f$ is ${\cal L}$-integrable on lots of subsets of this interval. On these sets $f$ can be
integrated twice  to give  contributions to the value of 
${\cal T}_{{(2,s)}_0}\mkern-4mu-\mkern-8mu\int_{(-\pi,\pi)}^x f$ that in the case of intervals are of the form (9).  The
difficulty is to put these contributions together correctly to get the value of the second order
primitive ${\cal T}_{{(2,s)}_0}\mkern-4mu-\mkern-8mu\int_{(-\pi,\pi)}^x f$.

The totalization process used here is  much more difficult than  the  totalization used for  classical primitive problem. Much
deeper knowledge of the fine properties of sets is required, as well as some very deep properties of the Schwarz
derivative. As a further evidence of the extra complications Denjoy distinguishes nine distinct types of limit calculations that  are needed to
complete the calculation, each used infinitely often in general--- and he gives examples to show that all are needed; for the
classical primitive problem he only distinguishes three, and all are very simple.

 These difficulties were overcome by Denjoy who went on to show how to complete the calculation and 
gave  formul\ae\  for the coefficients, like those in (F), but now more
complicated since we are using a second order 
integral,  rather than a first order  integral; see (14) below.

 To see the order of these difficulties  consider the simplest situation: we have computed the total  on two adjoining
intervals $[a,b]$ and $[b,c]$ how do we compute the value  on $[a,c]$? 

In the case of the classical primitive problem this means
that if 
$$
\int_a^x\!g = \Gamma_1(x)-\Gamma_1(a), \;a\le x\le b,\quad{\rm and}\quad \int_b^x\!g =\Gamma_2(x)-\Gamma_2(b), \;b\le x\le c,
 $$
then, as is  easy to see,
$$
\int_a^x\!g = \Gamma(x)-\Gamma(a), a\le x\le c,
$$
where 
$$
\Gamma(x) = 
\cases
{\Gamma_1(x), & if $a\le x\le b$,\cr
\Gamma_2(x) +\Gamma_1(b)- \Gamma_2(b), & if $b\le x\le c$.\cr
}
$$
However the situation  for the $({\cal T}_{2,s})_o$-integral  is much more complicated;  \ref{D, vol. III,
pp.278--279}.

Let
$$
{\cal T}_{{(2,s)}_0}\mkern-4mu-\mkern-8mu\int_{(a,b)}^x g = (x-a)(x-b)\bigl[a,x,b; G_1\bigr], a\le x\le b; 
$$
 and
$$
{\cal T}_{{(2,s)}_0}\mkern-4mu-\mkern-8mu\int_{(b,c)}^x g = (x-b)(x-c)\bigl[b,x,c; G_2\bigr], b\le x\le c.
$$
There is an $G$ such that $D_s^2G = g$ on $[a,c]$, which differs by a linear function from $G_1$ on $[a,b]$, and by a
linear function from $G_2$ on $[b,c]$. Assume for simplicity that $a\le x<b$, (the case $b<x\le x$
can be treated similarly);   then by (17) below
$$
\eqalign
{
{\cal T}_{{(2,s)}_0}\mkern-4mu-\mkern-8mu\int_{(a,c)}^x g =& (x-a)(x-c)\bigl[a,x,c; G\bigr]\cr
 =&  (b-c)(x-a)\bigl[a,b,c; G\bigr] +(x-a)(x-b)\bigl[a,x,b; G\bigr]\cr
=&  (b-c)(x-a)\bigl[a,b,c; G\bigr] +(x-a)(x-b)\bigl[a,x,b;G_1\bigr]\cr
=&  (b-c)(x-a)\bigl[a,b,c; G\bigr] +{\cal T}_{2,s}\mkern-4mu-\mkern-8mu\int_{(a,b)}^x g.\cr
}
$$
It remains to see how we can compute the second divided difference in first term in the last line, knowing only  $G_1$ and
$G_2$. 

Let $a<b-h<b<b+h<c$ and then
$$
\eqalignno
{
\bigl[a&,b,c; G\bigr]&(10)\cr
 =&(c-a)\left((b-h-a)\bigl[a,b-h,b; G\bigr] +(c-b+h)\bigl[b,b+h,c; G\bigr]  + \frac{\Delta^2
G(x;h)}{h}\right)\cr
=& (c-a)\left((b-h-a)\bigl[a,b-h,b; G_1\bigr] +(c-b+h)\bigl[b,b+h,c; G_2\bigr]  + \frac{\Delta^2
G(x;h)}{h}\right)\cr
}
$$
 Since $D_s^2G$ is finite $G$ is {\tt smooth}, and so  the last term on the right tends to zero with $h$. Hence we get that: 
 $$
\bigl[a,b,c; G\bigr]=(c-a)\lim_{h\to 0}\biggl((b-h-a)\bigl[a,b-h,b; G_1\bigr] +(c-b+h)\bigl[b,b+h,c;
G_2\bigr]\biggr);\eqno(11)
$$
of course this gives the value of ${\cal T}_{{(2,s)}_0}\mkern-4mu-\mkern-8mu\int_{(a,c)}^b f $.
 
This is a difficulty occurs  at the very first step of the calculations and there are many more
difficulties along the road.

\subsn{ The James $P^2$-integral}  Given the simplicity of the
Perron approach to the classical primitive problem, see {[Bu, 6.2]}, it would seem natural to ask whether  a
similar idea could not be used for  the solution of the of the coefficient problem.

 However it was not until the final appearance of Denjoy's work  in the fifties that 
R.D.James, on the suggestion of A. Zygmund, used the Perron method  on
the Schwarz derivative;  see \ref{J1;
 J-G; Z, vol. II, p.86}.
  
The extension of the Perron method to this situation is not immediate. The derivative being used
is a second order derivative and the corresponding integral, like Denjoy's in section 7.1, is a
second order integral

The classes of major and minor  functions defined in  \ref{Bu2, 6.2}  are generalized in  a natural way.  Call $M,m$
major and minor
functions for
$f$ on
$[a,b]$ if
$$
M(a) = m(a) = M(b) = m(b)=0\qquad \underline D_s^2M(x)\ge f(x)\ge \overline
D_s^2m(x),\; a<x<b;
$$
 and then the {\sl symmetric second  Perron,  or (James) ${\cal P}^2$-, integral of $f$\/ } is 
$$
{\cal P}^2\mkern-4mu-\mkern-8mu\int_{(a,b)}^xf= \sup M(x) = \inf m(x), a\le x\le b\eqno(12)
$$
 when these two families of functions have a common sup and inf.

 Any function, $F$, that differs from the
 right-hand side of (12) by a linear function will be called a {\sl (James) ${\cal P}^2$- integral of $f$\/}; and 
so by (8),
$$
{\cal P}^2\mkern-4mu-\mkern-8mu\int_{(a,b)}^xf=  (x-a)(x-b)\bigl[a,x,b; F\bigr]\eqno(13)
$$
 A  ${\cal P}^2$-integral $F$, of $f$,  can be proved to be  continuous, $ACG$, smooth,  differentiable almost everywhere, and   $D^2F=
f$ almost everywhere; see \ref{J1; S1}

However  a ${\cal P}^2$-integral need not be  differentiable everywhere; as the following example, [J1], shows:
$$
H(x) = \cases
{ x\cos \inv{x}, & if $x\ne 0$,\cr
0,& if $x=0$;\cr
}\qquad
h(x) = \cases
{
-x^{-3}\cos\inv{x}, & if $x\ne 0$,\cr
0,& if $x=0$.\cr
}
$$
 Then $H$ is continuous. smooth and $D_s^2 H= h$ and so  its  ${\cal P}^2$-integral is given by  (13),
with $H,h$ replacing  $f, F$ respectively. However $H'(0)$ does not exist, as is easily seen.

In addition the  ${\cal P}^2$-integral does not have nice properties with respect to additivity. That
is a function can be integrable both on $[a,b]$ and on $[b,c]$ but not on $[a,c]$. Consider the
following example due Skvortsov, \ref{Sk};
$$
g(x)= 
\cases
{
0, & if $ -2/ \pi\le x \le 0$,\cr
-x^{-3}\cos\inv{x}, & if $0<x\le 2/ \pi$.\cr 
}
$$
Obviously $g$ is integrable on $[-2/ \pi,0]$, with integral $0$, and by the previous
example it is integrable on  $[0, 2/\pi]$, with integral given by the function $H$ of that example.
Suppose then that it was integrable on the complete interval $[-2/ \pi,2/ \pi]$ with integral
$G$. Then $G$ would be linear on the left half of the interval, and differ from $H$ by a linear
function on the second half of the interval. This implies that $G$ has a left derivative at the
origin, but  that it does not have a right derivative at the origin, and so is not smooth there.  As
a result we deduce that
$g$ is not integrable on the complete interval.

Skvortsov gave a necessary and sufficient
condition for the additivity to hold; it is a little complicated to state; see [Sk] and [C2]. A 
simple   corollary is the sufficient condition: if $g$ is James integrable on $[a,b]$, and also on
$[b,c]$, with integrals $ G_1$ and $G_2$ respectively, then $g$ is integrable on  $[a,c]$ if $G_1$
has a left derivative at
$b$, and $G_2$ has a right derivative at $b$.

It is worth remarking that if the integral exists on $[a,b], [b,c]$ and on $[a,c]$ then the argument used in (10) and (11) can be
repeated to connect the values of the three integrals;  the function $F$ used there  is now the integral  $[a,c]$  which is smooth so
we can get (10) from (11) as was done for the ${\cal T}_{{(2,s)}_0}$- integral.

However the James integral does integrate second
symmetric derivatives, and so the sum 
  function of (T), $f$ say,  is integrable in this sense, and the sum of (7), $F$, is a James integral
of $f$.  Hence using the remarks in 11.1 about second divided differences
$$
a_0 = -\frac{1}{\pi^2}{\cal P}^2\mkern-4mu-\mkern-8mu\int_{(-2 \pi,2 \pi)}^0f.
$$
Further, using the method of formal multiplication,  each of $f(x)\sin kx, f(x)\cos kx,  k=1, 2,
\ldots$ are also integrable in this sense and
$$
\eqalign
{
a_k =&
-\frac{1}{\pi^2}{\cal P}^2\mkern-4mu-\mkern-8mu\intdl{(-2\pi,2\pi)}{0}{f(y)\cos ky}{y},\ k= 0,1,\ldots\, ,\cr
b_k =&
-\frac{1}{\pi^2}{\cal P}^2\mkern-4mu-\mkern-8mu\intdl{(-2\pi,2\pi)}{0}{f(y)\sin ky}{y},\ k= 1,\ldots\, .\cr
}
\eqno(14)
$$
(A comment on the signs in these formul\ae\ might be in order. If $f>0$ we would expect $a_0>0$, but then  $D^2F=f>0$, so $F$  would
be convex and so $\bigl[-2 \pi,0,2 \pi; F\bigr]>0$, see 11.1,  and $-4 \pi^2\bigl[-2 \pi,0,2 \pi;
F\bigr]={\cal P}^2\mkern-4mu-\mkern-8mu\int_{(-2 \pi,2 \pi)}^0f<0$.)
 
\sectn{Two Other  Integrals of Perron-type}

\asubsn{The Integral of Marcinkiewicz \&
Zygmund}
 The first trigonometric integral to solve the coefficient problem, after the work of
Denjoy was announced in the early twenties, was given by  Marcinkiewicz \&
Zygmund, \ref{M-Z}.  

They  took on the difficulties of the first primitive,
 getting  around the lack of continuity and the fact that $\Phi$, (1),  was defined only almost
everywhere in an ingenious way.  As we have seen $\Phi$ is Lebesgue  integrable, and to
integrate it we do not need it to be defined everywhere. Using  a
symmetric derivative called the   {\tt symmetric Borel derivative} that is based on
smoothing the function
to be differentiated by first integrating; \ref{Th3 pp.15--16},   they were able to solve the modification of the
classical primitive problem for this derivative and to show that the integral obtained, the so
called {\sl${\cal  T}$-integral\/},  solved the coefficient problem.  The details of the construction of the
${\cal T}$-integral  are very similar to those of the better known   integral due to Burkill to be described below.
Formal multiplication was used to obtain (F).

 \asubsn{ The Burkill  Integral}
 A trigonometric integral that is related to the James integral   was  defined
by J.C. Burkill, see \ref{BJ}. His approach has the advantages of that of the previous section as it results in a first order integral, but the tools used
are standard results from the theory  of trigonometric series, which makes this integral  more approachable than the  ${\cal
T}$-integral.

Consider the
following:  the sum of the series (7), $F$,  is the indefinite integral of the  function $\Phi$ in (1)
hence
$$
\eqalign
{
\frac{F(x+h) +F(x-h)- 2F(x)}{h^2}= &\frac{1}{h^2}\intdl{x}{x+h}{\Phi(y)}{y}
- \intdl{x-h}{x}{\Phi(y)}{y}\cr
=&\frac{1}{h^2}\intdl{0}{h}{\bigl(\Phi(x+u) - \Phi(x-u)\bigr)}{u}.\cr
}
$$
 So that 
$$
D_s^2 F(x) = \lim_{h\to 0}\frac{1}{h^2}\intdl{0}{h}{\bigl(\Phi(x+u) -
\Phi(x-u)\bigr)}{u}.\eqno(15) $$
 The  right-hand side of (15)  can be considered as firstly, by integrating, smoothing out the bad
function $\Phi$, and then differentiating, giving a generalized  first order derivative of 
$\Phi$,  called the {\sl symmetric Ces\`aro derivative of $\Phi$, at $x$\/}; see 11.3.

 Burkill's idea was to
use  this symmetric  Ces\`aro derivative   to define an  integral, using the
 Perron approach, and to show that it  solves the coefficient problem. In this way  the use the
second integrated series, (7), is avoided.

 As we saw above the first integrated
series does not necessarily converge everywhere and so, like the ${\cal T}$-integral,  the
integral defined by
Burkill, the {\sl symmetric  \ Ces\`aro -Perron\/}, or ${\cal SCP}$-, {\sl integral\/} has a primitive
that is only defined almost everywhere. This makes the  details of the construction, the
modifications needed to the definitions of the major and minor functions in this case,   very intricate and so will not be given here.

 Using this integral the formul\ae\ (F) were obtained by  an integration by parts formula proved by  Burkill.
The proof of this was found, by Skvorcov,  to contain a flaw that was only corrected much later,
by Sklyarenko; \ref{S3; C-Th }. In the meantime the formul\ae\ (F) were reproved using the
method of formal multiplication; \ref{BH}

If the correct formulation of the definition of the major and minor functions is made a function is ${\cal SCP}$-integrable if and
only if it is ${\cal P}^2$- integrable, but the discussion of this point is very technical; see \ref{V, p. 681}.

\sectn{The Story Continues}

As we have seen in the ${\cal T}$- and the  ${\cal SCP}$-integrals the difficulties with $\Phi$ are avoided by the
use of  symmetric  differentiation and  by using integration to smooth out bad functions.
However there is a much more direct modification of the classical primitive problem that can
be used.

Associated with $\Phi$ there is a symmetric derivative that plays a role analogous to that
played by the Schwarz  derivative for the function $F$.  

The   first order symmetric derivative of $\Phi$,  see (19), is not quite right since although it does not need  $\Phi$ to be defined
everywhere it does need it to be defined in
some interval around each point; and this  may not be the case. This difficulty is avoided by generalizing 
the derivative in (19) even further.

The limit process used in defining continuity and derivatives can be modified in many ways. The
most useful modification in this connection is the so-called {\sl approximate limit}. In this we
only require that 
$ x+h\to x$ through a set that has density 1 at $x$.  That is $x+h$ is restricted to a set $E_x$
with
$$
 \lim_{h\to 0}\frac{\big|E_x\cap[x-h, x+h]\big|}{2h} = 1,
$$
(where by $ |A|$ we mean the measure of $A$.) The theory of such limits, approximate   limits, 
can be developed and gives these limits all the properties of usual limits; see for instance [G, p.223; Z I, p.23]. The correct derivative for use here
is one in which the limit in (19) is taken in this approximate sense giving the {\sl approximate symmetric derivative of $\Phi$ at $x$}; 
$$
 D_{s,ap}\Phi (x) = {\rm ap}\!-\!\lim_{h\to 0} \frac{\Phi(x+h)-\Phi(x-h)}{2h},\eqno(16)
$$
where of course ${\rm ap}\!-\!\lim_{h\to 0} $ means the approximate limit as $h\to 0$. It is the
correct
derivative because of  the following result of Rajchman \& Zygmund,   [Th3, pp.17--18; Z, vol.I,  p.324]. 
\proclaimsc{}{If $f$ is the sum
of  the series  (T) and if $\Phi$ is as in (1) then
 $$ 
 D_{s,ap}\Phi  = f.
$$
}\par
This result was proved in 1926 and it is reasonable to ask why it could not
have been  used to give an an exceedingly simple generalization of the Perron integral, one with
major and minor functions defined  by
$$
M(a) = m(a) = 0\qquad\underline D_{s,ap}M\ge g\ge\overline D_{s,ap}m
$$
to get a trigonometric integral that would solve the coefficient problem?

A critical part of the Perron theory is the fact that for each major function $M$, and each minor function	 $m$, the difference $
M-m$ is an increasing function. 

In the
classical
Perron case this derives from a fairly elementary result: 
\proclaimsc{}{if $\underline Dh\ge 0$ then $h$ is increasing.}\par
This implies that in the classical case $\underline D(M-m)\ge 0$, and so   $
M-m$ is increasing.

The same result for the approximate symmetric 
derivative turns out to be much more difficult to prove,  was not known until very recently when Freiling \& Rinnie, \ref{F-R}, showed that:
 \proclaimsc{}{if $\underline
D_{s,ap}h\ge 0$ then $h$ is increasing.}\par

This allowed
Preiss \& Thomson to give another solution to the coefficient problem along the lines
suggested above; a Perron integral that inverts the approximate symmetric derivative; [P-Th2].

At the same time they also gave an equivalent Riemann definition of an integral based on the
Kurzweil-Henstock methods mentioned in \ref{Bu2}.. 

To do this they had to consider symmetric partitions of an interval;
that is partitions  in which we would want that $x_{i-1} = y_i -h_i$, $x_i= y_i
+h_i, h_i<\delta(y_i)$; this extension is far from trivial and is due to Thomson; [Th2], [P-Th1]. In
addition they had to consider approximate partitions that would have  the $x_i $  lie
in some set $E_{y_i}$  of density 1 at
$y_i$.  Such partitions were first considered by Henstock, but in the case of symmetric
partitions they are much more difficult to handle; [P-Th2]. The details as is always the case
with the coefficient problem are complicated. Finally they also gave a variational form for this
integral. A very detailed discussion of all these integrals can be found in \ref{Th3}.

\sectn{Final Remarks}

 Firstly we note that all the integrals introduced extend the ${\cal D^*}$-integral although  relationships with other general classical integral is
more complicated; see \ref{V}.

 Do these solutions to the  classical coefficient  problem  end of the story ? No,  there are at several  other areas of further
interest.

An obvious question is,  are all the solutions above equivalent, as the various
solutions to the  classical primitive problems are. In general the answer to
this is unknown. Using  the classical definitions the  ${\cal P}^2$-integral is more general than the  ${\cal SCP}$-integral but    as we have
pointed out  these  integrals are essentially equivalent if some changes are made in the definitions;  although even here there are technical problems.
A trigonometric integral that lies  between the  ${\cal SCP}$-integral and the ${\cal P}^2$-integral has recently been given by Mukhopadhyay,
\ref{M}.   Just whether any of the various integrals are equivalent to the original
${\cal T}_{{(2,s)}_0}$-integral is unknown.
 A very detailed examination of this problem can be found in the recent masterly paper by
Cross \& Thomson, [C-Th]

Then there is no end to the class of integrals that can be
investigated by considering other methods of convergence for (T). Equality
in (T) has so far meant everywhere convergence of the series to the function, but
we could instead consider some  summability method, as is very common
in the theory of trigonometric  series. The use of Ces\`aro  summability leads to a 
family of integrals that generalize the James integral described above, [J2--4];
also, less successfully the Burkill  integral has been extended to a class  of integrals that does  a
similar job,  [Bu-L].   A mention should also be made of an
attempt to give a Riemann form to these higher order integrals, [C2]. A full discussion of all of  these
integrals has not been given and the technical problems seem to be daunting.

In addition an integral to be used when Abel summability is
considered has been defined by Taylor; this integral is more general than the ${\cal  P}^2$-integral; see [C1; T].

By far the best reference for this material is the excellent  book by Brian S.Thomson, [Th3], where
all the facts are sorted out and given a very lucid exposition. However it not a book for the light
-hearted. The article by James, [J3], is an easy introduction to the field although a little dated;
the same can be said of the lecture by R.L.Jeffery, [Je]. A short but excellent discussion can also
be found in the  article by Vinogradova \& Skvorcov, [V],  as well as in the book of Zygmund, where it is preceded by an equally
short review of the classical primitive problem; [Z II, pp.83--90].
\sectn{ Appendix}

\bsubsn{Divided Differences.} 
 The {\sl first divided difference of $G$ at $a,x$\/} is written $[a,x;G]$,  where 
$$
[a,x;G]=\frac{
G(x) -G(a)}{x-a};
$$
The {\sl second divided difference of $G$ at\/} $a,x,b$ $\bigl[a,x,b; G\bigr]$ is  written $[a,x,b; G]$, where
$$
\eqalign
{
\bigl[a,x,b; G\bigr] =& \frac{1}{b-a}\biggl\{\frac{G(b)-G(x)}{b-x}\,- \,\frac{G(x)-G(a)}{x-a}\biggr\}\cr
=&\frac{G(a)}{(a-x)(a-b)}+\frac{G(x)}{(x-b)(x-a)}+\frac{G(b)}{(b-a)(b-x)}.\cr
}
$$

 A few very simple facts about the second divided difference are worth noting.

(i)  if $h$ is a polynomial of  degree at most two then for all $a,x,b,\; \bigl[a,x,b; h\bigr] =\half
h''$.

 In particular if $h(x) = Ax^2, \;  \bigl[a,x,b; h\bigr] =A$.

(ii) If $h(x) = \sin nx$ or $\cos nx,\; n= 0, 1,2,\ldots $ then $\bigl[ -2\pi,0,2\pi; h\bigr] =
0$.

(iii) If $a,b,c,x$ are any four distinct points then
$$
(c-x)\bigl[a,x,c; G\bigr]=(c-b )\bigl[a,b,c; G\bigr]  + (b-x)\bigl[a,x,b; G\bigr]. \eqno(17)
$$
(iv) A function $G$ is convex on [$a,b]$ if and only if $\bigl[x,y,z; G\bigr]\ge 0$ for all distinct $x,y,z$ in
$[a,b]$.

\bsubsn{ Odd and Even Parts of  Function} 
Any function $G$ can be written  in terms of its {\sl odd and even parts at \/} $x$
$$
G(x+t)= \frac{G(x+t) + G(x-t)}{2} + \frac{ G(x+t) -G(x-t)}{2} = G_e^x(t) + G_o^x(t)
$$
 and the derivatives of  these two parts at $t=0$ depend on
$$
\eqalign
{
 G_o^x(t)-G_o^x(0) =& \frac{G(x+t) -G(x-t)}{2};\cr
G_e^x(t)-G_e^x(0) =&\frac{G(x+t) + G(x-t)-2G(x)}{2}= \half\Delta^2f(x;t).\cr
}
\eqno(18)
$$
So the derivative at the origin of the odd part of a function  exists precisely when the function has a  {\sl first order
symmetric derivative at \/}  $x$:
 $$
D_sG(x) = \lim_{h\to 0}\frac{G(x+h) - G(x-h)}{2h}.\eqno(19)
$$
This derivative certainly exists when $G'(x)$ exists but not conversely as taking $G(x) = |x|$ shows; more, the value of
$G$ at $x$ is not used in  (19) so  $D_sG(x)$  can exist even if $G$ is not continuous at $x$, or even if it not defined
there; consider for instance $ G(x) = \cos 1/x, x\ne 0$, when $D_sG(0)= 0$. See \ref{Th3, pp.5--7}.

 The second of the quantities  in (18)  was introduced in \ref{G IV, p.17}  when
defining the {\sl  Schwarz\/}, or {\sl second  symmetric derivative of a function $G$ at\/} $x$, 
$$
 D_s^2GF(x) = \lim_{h\to 0}\frac{G(x+h) + G(x-h) -2G(x)}{h^2},\eqno(20)
$$
 The derivative at the origin of the even part of a function  exists precisely when the function has the property
$$
 \lim_{h\to 0}\frac{G(x+h) + G(x-h) -2G(x)}{h}= 0; 
$$
 the function   $G$ is then said to be {\sl smooth at \/} $x$.
If $G$ is smooth at $x$ and if $G_+'(x)$  exists  then so does $G_-'(x)$ and $G_-'(x)=G_+'(x)$; of course a function can
be smooth without being continuous, let alone differentiable. The concept goes back to Riemann and was given the name by Zygmund;
\ref{Th3, p.159; Z, vol.I, p. 43}. 
     
\bsubsn{ Derivatives Using Integrals}  The idea of using integrals to smooth a function before taking the
derivative goes back to Borel, \ref{Th3, pp.15--16}: his {\sl mean value derivative\/}  is
$$
D_Bf(x)=\lim_{h\to 0}\recip{h}\intdl{0}{h}{\frac{f(x+t)-f(x)}{t}}{t},
$$
and the {\sl symmetric Borel derivative\/} , the one used for the ${\cal  T}$ -integral in 8.1, is
$$
D_{s,B}f(x)=\lim_{h\to 0}\recip{h}\intdl{0}{h}{\frac{f(x+t)-f(x-t)}{2t}}{t},
$$
Burkill introduced his  Ces\`aro derivative,
$$
CDf(x)=\lim_{h\to 0}\frac{2}{h^2}\intdl{0}{h}{\left(f(x+t)-f(x)\right)}{t}.
$$
and used the symmetric version, (15), for the ${\cal SCP}$-integral.

Of course given the variety of integrals at our disposal it is natural to ask which integral should be used in these
definitions. It turns out to make some difference. Burkill used the very natural  ${\cal D}^*$-integral, but by using
the
${\cal DH}$-integral instead Sklyarenko was able to show that the more general ${\cal SCP}$-integral thus obtained 
was equivalent to the ${\cal P}^2$-integral; \ref{S1}.
 \sectn{References}

\napa{[Bu]\ P S Bullen} The search for the primitive,  {\sl submitted\/}.
\napa{[Bu-L]\ P S Bullen \& C.M.Lee} The $SC_nP$-integral and
the $P^{n+1}$-integral, {\sl\ Canad.\ J.\ Math.\/}, 25\ (1973),\ 1274
--1284.
\napa{[BH]\ H  Burkill }  Fourier series of $SCP$-integrable
functions,{\sl\ J.\ Math.\ Anal.\ Appl.\/}, 57\ (1977),\
587--609.
\napa{[BJ]\ J C Burkill} Integrals and trigonometric series, {\sl\ Proc.\ London Math.\
Soc.},\enspace  (3)\ 1\ (1951),\ 46--57. Corrigendum: {\sl\ Proc.\ London Math.\ Soc.\/},
(3)\  47\ (1983),\ 192.
\napa{[C1]\ G E Cross}  On the generality of the $AP$-integral, {\sl\ Canad.\ J.\ Math.\/},
23\ (1971),\ 557-561.
\napa{[C2]\ G E Cross} Generalized integrals as limits of Riemann-like
sums, {\sl\  Real Anal.\ Exchange\/}, 13\ (1987--1988),\
390--403.
\napa{[C-Th]\ G E Cross \& B S Thomson} Symmetric integrals and trigonometric series,
{\sl Diss. Math,\/},\enspace CCCXIX,\ 1992.
\napa{[D]\ A Denjoy}{\it Le\c cons sur le Calcul de Coefficients d'une S\'erie
Trigonom\'etrique, I--IVa,b\/}, \ Gauthier-Villars,\  Paris 1941, 1949.  
\napa{[F-R]\ C Freiling \& D Rinnie} A symmetric density property; monotonicity and the
approximate
symmetric derivative, {\sl\  Proc.\ Amer.\ Math.\ Soc.\/}, 104\ (1988),\ 1098~--~1102.
\napa{[G]\ A  D  Gluchoff} Trigonometric series and theories of integration, {\sl\
Math.Mag.\/}, 67\ (1994),\ 3--20.
\napa{[Go]\ R A Gordon}  {\it The Integrals of Lebesgue, Denjoy, Perron and Henstock\/},  Amer.\
Math.\ Soc Memoir, 1994. 
\napa{[H]\ G  H  Hardy}{\it Divergent Series\/},  Oxford University Press, 1949. 
\napa{[Ho]\ E W Hobson}{\it The Theory of Functions of a Real Variable 
and the Theory of Fourier's Series I, II\/}, Cambridge University Press, 1926 
\napa{[J1]\ R D James} A generalized integral
II, {\sl\ Canad.\ J.\ Math.\/}, 2\ (1950),\ 297--306. 
 \napa{[J2]\ R D James} Generalized $nth$\
primitives, {\sl\ Trans.\ Amer.\ Math.\ Soc.\/},76\ (1954),\ 149--176.
 \napa{[J3]\ R D James} Integrals and
summable trigonometric series, {\sl\ Bull.\ Amer.\ Math.\ Soc.\/, 61 (1955),\
1--15.
\napa{[J4]\ R D James}
Summable trigonometric series {\sl\ Pacific\ J.\ Math.\/}, 6\ (1956),\ 99--110.
\napa{[J-G]\ R D James \&\ W H Gage} A generalized integral,{\sl\ Trans.\ Roy.\
Soc.\ Canada\/},\enspace III, (3)\ 40\ (1946),\ 25--35.
\napa{[Je]\ R L Jeffery}  Trigonometric Series, {\sl\ Canad.\ Math.\ Congress Lecture
Ser.\/}, \#~2,\ Toronto,\ 1956.
\napa{[K]\ K Knopp}{\it Theory and Application of Infinite Series\/}, Blackie \& Son Ltd., \ London,\
1948. 
\napa{[L]\ H  Lebesgue}{\it  Le\c cons sur l'Int\'egration et la Recherche des
Fonctions Primitives\/},  Gauthiers--Villars,  Paris; Ist Ed. 1904; 2nd Ed.1928.
\napa{[M-Z]\ J Marcinkiewicz \& A Zygmund} On the differentiability		of functions
and the summability of series, {\sl\  Fund.\ Math.},\enspace  26\ (1936),\  1--43.
\napa{[M] S N Mukhopadhyay} An extension of the $SCP$-integral with a relaxed integration by parts formula,
{\sl Analysis Math.\/}, \enspace 25\ (1999),\ 103--132.
\napa{[N]\ P Nalli}  Sulle serie di Fourier delle funzione non assolutamente integrabili,{\sl\ Rend.\ Circ.\
Math.\ Palermo\/}, 40\ (1915),\ 33--37.
New York, 1964
\napa{[P] I N Pesin}{\it Razvite Ponyatiya Integrala\/}\noter{\sevencyr  I N Pesin, Razvite Ponyatiya Integrala.} ,
Moscow,\enspace 1966. {\sl Engl.\
 transl.}\enspace : {\it Classical and Modern Integration Theory\/},\noter{\sevenrm Care using this  translation, several  terms are used in
their historical,  rather than their modern senses.} New York, 1970.
\napa{[P-Th1]\ D Preiss \& B S Thomson} A symmetric covering
theorem, {\sl\  Real Anal.\ Exchange\/},14\ (1988--1989),\
253--254.
\napa{[P-Th2]\ D Preiss \& B S Thomson}  An approximate
symmetric integral, {\sl\ Canad.\ J.\ Math.\/}, 41\ (1989),\
508-555.
 \napa{[S1]\ V A Sklyarenko}\kern-3pt\noter{\cyrf V.A.Sklyarenko}  Nekotorie svo\u \i stva $P^2$-primitivno\u \i,{\sl\ Mat.\
Zametki\/}, 12\ (1972),\ 693--700. {\sl Eng.\ transl.:\ Math.\ Notes\/}, 12\ (1972),\ 856--860,\
(1973).
 \napa{[S2]\ V A Sklyarenko} On Denjoy integrable
sums of everywhere convergent trigonometric	 series, {\sl\ Dokl.\ Akad.\ Nauk SSSR\/},
210\ (1973),\ 533--536. {\sl\ Engl.\ transl.:Soviet Math.\
Dokl.\/}, 14\ (1973),\ 771--775.
\napa{[S3]\ V A Sklyarenko} Ob integrirovani\u\i\ po \v castyam v $SCP$-integrale Burkillya, {\sl\
Mat.Sb.}, (154)\ 112 ,\ (1980),\ 630--646.{\sl\ Engl.\ transl.\ Math.\ USSR
Sbornik}, 40\ (1981),\ 567--583.
\napa{[S4]\ V A Sklyarenko} On a property of the Burkill SCP-integral, {\sl  Mat.\ Zamet.\/}, 65(1999),\ 599-604. 
\napa{[Sk]\ V A Skvorstov}\kern-5pt\noter{\cyrf V.A.Skvortsov. {\sevenrm Also transliterated as Skvorcov.}} Po povody opredeleni\u\i\
$P^2$- i $SCP$-integralov, {\sl\ Vestnik Moscov.\ Univ.; Ser.I Mat.\ Mehn.\/}, 21\ (1966),\ 12--19.
\napa{[T]\ S J Taylor} An integral of Perron's type defined with the help of trigonometric
series, {\sl\  Quart.\ J.\ Math.\ Oxford\/}, (2)\ 6 \ (1955),\ 255--274.
\napa{[Th1]\ B S Thomson}  Symmetric derivatives and symmetric
integrals, {\sl\  Real Anal.\ Exchange\/}, 15\ (1989--1990),\
49--61.
\napa{[Th2]\ B S  Thomson}Some symmetric covering lemmas,
{\sl\  Real Anal.\ Exchange\/}, 15\ (1989--1990),\ 346--383.
\napa{[Th3]\ B S  Thomson}{\it Symmetry Properties of Real Functions\/}, Marcel
Dekker Inc.,~1994 .
\napa{[V]\ I A Vinogradova \& V A Skvorcov}\kern-5pt\noter{\cyrf I.A.Vinogradova} 
Generalized Fourier integrals and series,{\sl\ J.Soviet Math.\/}, 1\ (1973),\ 677--703.
\napa{[Y]\ W H Young} On non-absolutely
convergent, not necessarily continuous, integrals, {\sl\ Proc.\ London Math.\ Soc.\/}, (2)\
16 \ (1917),\ 175--218.
\napa{[Z]\ A Zygmund}{\it Trigonometric Series I, II\/}, Cambridge,\ 1959.

\bye